\documentclass[12pt]{amsart}
\usepackage{amssymb}

\begin{document}
\oddsidemargin=-1cm
\evensidemargin=-1cm
%\hoffset 5.5truemm
%\setlength{\textwidth}{25cm}

\title[Dynamics of Segre varieties]{Dynamics of the
Segre varieties of a real submanifold in complex
space}
\author[M. S. Baouendi, P. Ebenfelt, L. P.
Rothschild]{M. S. Baouendi, P. Ebenfelt, and Linda
Preiss Rothschild}
\address{Department of Mathematics, 0112, University
of California at San Diego, La Jolla, CA 92093-0112}
\email{sbaouendi@ucsd.edu, lrothschild@ucsd.edu,
\hfil\break ebenfelt@math.kth.se}
%\abstract \endabstract
%\keywords \endkeywords
%\subjclass{32H02}
\thanks{2000 {{\em Mathematics Subject Classification.}
32H02}}
%\date{\number\year-\number\month-\number\day}
%\enddate
%\loadeufm

%\def\Label#1{\label{#1}{\bf \hbox{  }(#1)\hbox{  } }}
\def\Label#1{\label{#1}}
\def\1#1{\ov{#1}}
\def\2#1{\widetilde{#1}}
\def\3#1{\mathcal{#1}}
\def\4#1{\widehat{#1}}

\def\s{s}
\def\k{\kappa}
\def\ov{\overline}
\def\span{\text{\rm span}}
\def\ad{\text{\rm ad }}
\def\tr{\text{\rm tr}}
\def\xo {{x_0}}
\def\Rk{\text{\rm Rk\,}}
\def\sg{\sigma}
\def \emxy{E_{(M,M')}(X,Y)}
\def \semxy{\scrE_{(M,M')}(X,Y)}
\def \jkxy {J^k(X,Y)}
\def \gkxy {G^k(X,Y)}
\def \exy {E(X,Y)}
\def \sexy{\scrE(X,Y)}
\def \hn {holomorphically nondegenerate}
\def\hyp{hypersurface}
\def\prt#1{{\partial \over\partial #1}}
\def\det{{\text{\rm det}}}
\def\wob{{w\over B(z)}}
\def\co{\chi_1}
\def\po{p_0}
\def\fb {\bar f}
\def\gb {\bar g}
\def\Fb {\ov F}
\def\Gb {\ov G}
\def\Hb {\ov H}
\def\zb {\bar z}
\def\wb {\bar w}
\def \qb {\bar Q}
\def \t {\tau}
\def\z{\chi}
\def\w{\tau}
\def\Z{\zeta}

\def \T {\theta}
\def \Th {\Theta}
\def \L {\Lambda}
\def\b {\beta}
\def\a {\alpha}
\def\o {\omega}
\def\l {\lambda}

\def \im{\text{\rm Im }}
\def \re{\text{\rm Re }}
\def \Char{\text{\rm Char }}
\def \supp{\text{\rm supp }}
\def \codim{\text{\rm codim }}
\def \Ht{\text{\rm ht }}
\def \Dt{\text{\rm dt }}
\def \hO{\widehat{\mathcal O}}
\def \cl{\text{\rm cl }}
\def \bR{\mathbb R}
\def \bC{\mathbb C}
\def \C{\mathbb C}
\def \bL{\mathbb L}
\def \bZ{\mathbb Z}
\def \bN{\mathbb N}
\def \scrF{\mathcal F}
\def \scrK{\mathcal K}
\def \mc #1 {\mathcal {#1}}
\def \scrM{\mathcal M}
\def \cR{\mathcal R}
\def \scrJ{\mathcal J}
\def \scrA{\mathcal A}
\def \scrO{\mathcal O}
\def \scrV{\mathcal V}
\def \scrL{\mathcal L}
\def \scrE{\mathcal E}
\def \hol{\text{\rm hol}}
\def \aut{\text{\rm aut}}
\def \Aut{\text{\rm Aut}}
\def \J{\text{\rm Jac}}
\def\jet#1#2{J^{#1}_{#2}}
\def\gp#1{G^{#1}}
\def\gpo{\gp {2k_0}_0}
\def\emmp {\scrF(M,p;M',p')}
\def\rk{\text{\rm rk\,}}
\def\Orb{\text{\rm Orb\,}}
\def\Exp{\text{\rm Exp\,}}
\def\Span{\text{\rm span\,}}
\def\d{\partial}
\def\D{\3J}
\def\pr{{\rm pr}}

\def\dbl{[\hskip -1pt [}
\def\dbr{]\hskip -1pt]}

\def \D{\text{\rm Der}\,}
\def \Rk{\text{\rm Rk}\,}
\def \ima{\text{\rm im}\,}

\newtheorem{Thm}{Theorem}[section]
\newtheorem{Def}[Thm]{Definition}
\newtheorem{Cor}[Thm]{Corollary}
\newtheorem{Pro}[Thm]{Proposition}
\newtheorem{Lem}[Thm]{Lemma}
\newtheorem{Rem}[Thm]{Remark}

\maketitle

\section{Introduction}

Let $M\subset\bC^N$ be a smooth ($C^\infty$) real
submanifold of codimension
$d$ with $0\in M$. We choose smooth real-valued
functions
$r=(r_1,\ldots, r_d)$, with 
differentials $dr_1,\ldots, dr_d$  linearly 
independent at $0$, so that $M$ is defined by $r=0$
near the origin. If the complex differentials
$\partial r_1,\ldots,
\partial r_d$ are also linearly independent at $0$,
then
$M$ is called {\it generic} (near the origin). If $M$
is a real-analytic, generic submanifold, there is a
family of complex submanifolds of $\bC^N$, called
the Segre varieties associated to $M$, which carry a
great deal of information about the local geometry
of $M$. The Segre varieties have been used
 by many mathematicians to study
mappings between generic submanifolds. (See the end
of this introduction for some specific references.)
 In this paper, we shall consider an algebraic
substitute for these varieties for the case of smooth
manifolds by introducing a formal mapping which, in
the real-analytic case, parametrizes the Segre
varieties.  Our main objective
is to study iterations of this mapping and relate
these iterations to the local CR geometry of the
manifold. 

If $r$ is a local defining function of a smooth
generic submanifold $M$ as above, then we denote by
$\rho_j(Z,\bar Z)$ the Taylor series of $r_j$ at
$0$. We write $\rho=(\rho_1,\ldots,\rho_d)$. We
consider
$\rho_j(Z,\zeta)$ as a formal power series in the
$2N$ indeterminates
$(Z,\zeta)$. We shall denote the ring of such power
series with complex coefficients by
$\bC\dbl Z, \Z \dbr$. By a formal mapping
$F\colon(\bC^k,0)\to(\bC^p,0)$, we shall mean a
$p$-tuple
$F(x)=(F_1(x),\ldots,F_p(x))$, where $x=(x_1,\ldots,
x_k)$, of formal power series
$F_j\in\bC\dbl x\dbr$ without constant terms. The
rank of $F$, $\Rk F$, is defined as the rank of the
Jacobian  matrix $\d F /\d x$ regarded as a
$\Bbb K_x$-linear mapping
$\Bbb K_x^k\to\Bbb K_x^p$, where $\Bbb K_x$ denotes
the field of fractions of $\bC\dbl x\dbr$. Hence
$\Rk F$ is the largest integer $s$ such that there
is an
$s\times s$ minor of the matrix $\d F/\d x$ which is
not  0 as a formal power series in $x$.

Let $\gamma(\Z, t)$, where
$\Z=(\Z_1,\ldots\Z_N)$, $t=(t_1,\ldots,t_n)$, and
$n=N-d$, be a formal mapping
$(\bC^N\times\bC^n,0)\to (\bC^N,0)$ such that
\begin{equation}\Label{gamma}
\rho(\gamma(\zeta,t),\zeta)= 0,\quad
\rk\, \frac{\d\gamma}{\d t}(0,0)=n.
\end{equation} The existence of such $\gamma(\Z,t)$
is a consequence of the formal implicit function
theorem and the fact that
$\partial_Z\rho_1,\ldots\partial_Z
\rho_d$ are linearly independent at $0$. We shall
call a formal mapping $\gamma(\zeta,t)$ satisfying
(\ref{gamma}) a {\it Segre variety mapping} for the
germ of $M$ at $0$.  Recall that if $M$ is
real-analytic, in which case we may assume that the
formal series $\rho_j(Z,\zeta)$ are convergent, then
the Segre variety of $M$ at $p$, for $p$ near $0$,
is the complex
$n$-dimensional submanifold defined by the equation
$\rho(Z,\bar p)=0$.  Hence, in this case,
$\gamma(\zeta, t)$ satisfying (\ref{gamma}) can be
chosen to be convergent, and the mapping $t\mapsto
\gamma(\zeta,t)$, for $t$ near
$0\in\bC^n$,  parametrizes the Segre variety of $M$
at $\bar\Z$.

We define a sequence of formal mappings $v^j\colon
(\bC^{nj},0)\to (\bC^N,0)$, called {\it the iterated
Segre mappings} of
$M$ at $0$ (relative to $\gamma$), inductively as
follows:
\begin{equation}
\begin{aligned}\Label{vj} &v^1(t^1):=
\gamma(0,t^1),\\ &v^{j+1}(t^1,\ldots,t^{j+1}):= 
\gamma(\bar v^{j}(t^1,\ldots, t^j),t^{j+1}).
\end{aligned}
\end{equation}

Recall that $M$ is said to be of {\it finite type}
at $0$ (in the sense of Kohn \cite{Kohn} and
Bloom--Graham \cite{BG}) if the Lie algebra $\frak
g_M$ generated by the $(1,0)$ and the $(0,1)$ vector
fields tangent to
$M$ span $\bC T_0M$, the complexified tangent space
of $M$ at $0$ (see e.g.\ \cite {BER}, Chapter
I). The following is one of the main results of this
paper.

\begin{Thm}\Label{main} Let $M\subset \bC^N$ be a
smooth generic submanifold of codimension $d$ with
$0\in M$, and let
$\gamma(\zeta,t)$ be a Segre variety mapping of
$M$ at $0$. Let $v^j$, $j\geq 1$, be the iterated
Segre mappings of
$M$ at $0$ relative to $\gamma$. Then the rank $\Rk
v^j$ is an increasing function of $j$ and is
independent of the choice of the holomorphic
coordinates
$Z$, the defining function $r$, and the Segre variety
mapping $\gamma$.  In addition, there exists an
integer $k_0$,
$1\le k_0\leq d+1$, such that $\Rk v^j=\Rk v^{j+1}$
for
$j\ge k_0$, and if $k_0> 1$, then $\Rk v^j < \Rk
v^{j+1}$ for $1\le j\le k_0-1$.
 Moreover, the following are equivalent:
\begin{enumerate}
\item [(i)] $M$ is of finite type at $0$.
\item [(ii)]  $\Rk v^{k_0}=N$.
\end{enumerate}
\end{Thm}

For $p\in M$, we denote by $\frak g_M(p)$ the
subspace of
$\bC T_p M$ obtained by evaluating the vector fields
in
$\frak g_M$ at $p$.  We shall prove the following
result, which is more general than the equivalence
of (i) and (ii) in Theorem
\ref{main}.
\begin{Thm}\Label{main-2} Let $M\subset \bC^N$ be a
smooth generic submanifold of codimension $d$
through the origin, with
$\gamma$, $v^j$, and $k_0$ as in Theorem
$\ref{main}$. Let
$e=2N-d-\dim_\bC\frak g_M(0)$. Then, 
\begin{equation}\Label{main-id}
\Rk v^{k_0}=N-e.
\end{equation}
\end{Thm}

We also have the following result.

\begin{Thm}\Label{main-3} Let $M\subset \bC^N$ be a
smooth generic submanifold of codimension $d$
through the origin, with
$\gamma$, $v^j$, and $k_0$ as in Theorem
$\ref{main}$, and $e$ as in Theorem
$\ref{main-2}$.  Then there exist formal power series
$f_1,\ldots, f_e\in\bC\dbl Z\dbr$ such that
\begin{enumerate}
\item [(i)]
$df_1(0),\ldots,df_e(0),dr_1(0),\ldots,dr_d(0)$ are
linearly independent, where $r_1,\ldots,r_d$ are
defining functions for $M$ near $0$.
\item[(ii)] $f_k\circ v^j =0$, for $k=1,\ldots, e$
and
$j=1,2,\ldots$.
\end{enumerate}
\end{Thm}

We would like to point out that if $M$ is a
real-analytic generic submanifold, then one can
choose the Segre variety mapping $\gamma$ to be
holomorphic in a neighborhood of $0$. Hence, the
iterated Segre mappings $v^j$ are holomorphic. Also,
in this case, the power series $f_1,\ldots, f_e$ in
Theorem
\ref{main-3} may be chosen to be convergent near
$0$. The real submanifold in $M$ defined by
$f_1=\cdots=f_e=r_1=\cdots=r_d=0$ is the {\it local
CR orbit} of $0$ in $M$, and the complex submanifold
in $\bC^N$ defined by
$f_1=\cdots=f_e=0$ is its {\it intrinsic
complexification}. (The reader is referred e.g.\ to
\cite{BER}, Chapter I, for the definition of the
intrinsic complexification of a real-analytic CR
submanifold, and to \cite{BER}, Chapter IV, for the
definition of the local CR orbit.)

A smooth real submanifold $M\subset\bC^N$, defined
locally near
$p_0\in M$ by $r_1=\cdots=r_d=0$, is said to be
{\it CR} (near $p_0$) if the rank of the complex
differentials $\d r_1,\ldots,\d r_d$ is constant on
$M$ (near $p_0$). A real-analytic CR submanifold is
generic as a real submanifold of its intrinsic
complexification (see \cite{BER}, Corollary 1.8.10),
and hence Theorems
\ref{main}--\ref{main-3} also imply directly results
for real-analytic CR submanifolds. A merely smooth
CR submanifold need not be contained in any proper
complex submanifold, however, a formal version of
such an inclusion may be given in this case. (See \S 
\ref{Prelim} and Remark \ref{CR-ideal}.) Analogs of
Theorems \ref{main}--\ref{main-3} can be
formulated in this context, but we shall not do so
here.

As mentioned above, in the real-analytic case, the
use of Segre varieties for the study of
mappings of hypersurfaces and generic submanifolds has
a long history.  In particular we cite here the work
of Webster
\cite {Webster}, Diederich-Webster \cite {DW},
Baouendi-Jacobowitz-Treves \cite{BJT}, 
Diederich-Fornaess \cite{DF}, Baouendi-Roth\-schild
\cite{BR}, Forstneric~\cite{Forst}, Huang
\cite{Huang}.

The iterated Segre
mappings $v^j$ in the real-analytic case were
introduced, with a special choice of coordinates
(so-called canonical coordinates)
 and a special choice of Segre
variety mapping in \cite{Acta}. In that
paper, the image of the mapping
$v^j$ was called the
$j$th Segre set, and the existence of an integer
$k_0$ such that (\ref{main-id}) holds was
established. Also, for $M$ merely smooth, the
existence of the integer $k_0$ such that (i)
implies  (ii) in Theorem
\ref{main} can be deduced from  \cite {Acta} by
again using the special canonical coordinates
mentioned above (see \cite{MA} and \cite{JAMS}).
However, the results in Theorem
\ref{main}, even in the case where $M$ is
real-analytic, are sharper than those mentioned
above.  We should emphasize here that the proofs in
this paper are new, even in the real-analytic case,
and the approach here is more elementary than that in
\cite{Acta}; in particular, the proofs in this paper
do not assume the existence of canonical
coordinates. The results in Theorems
\ref{main}--\ref{main-3} above follow from a more
general statement, Theorem \ref{formal-main} below,
concerning formal manifolds. All the results are
thus reduced to questions about ideals and their
derivations in the ring of formal power series. 

The iterated Segre mappings, in the real-analytic and
algebraic cases, have played a crucial role in
recent work on mappings between generic
submanifolds. We mention here work of  
Zaitsev \cite{Z1}--\cite{Z2}, joint work of the first
and third authors with Zaitsev \cite{BRZ},
Mir \cite{Mirh}--\cite{Mirg}, as well as
work of the authors \cite{Acta}, \cite{CAG}, 
\cite{MA}, \cite{JAMS}. We should also mention the
work of Christ--Nagel--Stein--Wainger \cite{CNSW} in
which the ranks of iterated real mappings are related
to curvature and finite type conditions in a
very different context. Finally, we note that
in the real-analytic case another approach to the
construction of a sequence of mappings
$v^j$ for which  the identity (\ref{main-id}) holds
for some
$k_0$ was given by Merker \cite{Merker}.

The authors are grateful to Dmitri Zaitsev for showing us a
simple proof of Nagano's theorem, which was adapted to
the formal case in the proof of Proposition \ref{Nagano}.

\section{Preliminaries and reformulation of the main
results for  formal manifolds}\Label{Prelim}

In this section, we shall give some preliminary
material on formal manifolds and their mappings.  As
in the introduction, we let $\bC\dbl x\dbr$ denote
the ring of formal power series in $x=(x_1,\ldots,
x_k)$, and
$\Bbb K_x$ the field of fractions of $\bC\dbl
x\dbr$. An ideal $I\subset \bC\dbl x\dbr$ is called
a {\it manifold ideal} if, for some choice of
generators $f_1,\ldots, f_d$ of $I$, the vectors
$\d f_1/\d x (0),\ldots,
\d f_d/\d x (0)
$ are linearly independent in $\bC^k$.
(Observe that any other minimal set of generators
of $I$ must have $d$ elements whose differentials at
the origin are linearly independent.) If this
condition is satisfied we say that the manifold ideal
$I$ has codimension $d$ and dimension
$k-d$. (This terminology agrees with that of Krull
dimension for manifold ideals.) We shall say that
the  manifold ideal $I$ defines a {\it formal
manifold}
$\Sigma\subset\bC^k$ of dimension $k-d$ and
codimension $d$ and write
$I=\mathcal I(\Sigma)$. To motivate this
terminology, we observe that if the generators
$f_1,\ldots, f_d$ of
$I$ can be taken to be convergent, then the equations
$f_1(x)=\ldots=f_d(x)=0$ define a complex submanifold
$\Sigma$ through the origin in
$\bC^k$ (of codimension $d$) such that the ideal of
germs at $0$ of holomorphic functions vanishing on
$\Sigma$ is generated by $f_1,\ldots,f_d$. If
$\Sigma_1,\Sigma_2\subset\bC^k$ are formal manifolds
and
$\mathcal I(\Sigma_2)\subset\mathcal I(\Sigma_1)$,
then we shall say that $\Sigma_1$ is contained in
$\Sigma_2$ and write
$\Sigma_1\subset\Sigma_2$.

Recall that a formal mapping $F\colon (\bC^k,0)\to
(\bC^p,0)$ is a $p$-tuple $F=(F_1,\ldots, F_p)$ of
formal power series in $x$ that vanish at $0$. Any
formal mapping induces a $\bC$-linear
ring homomorphism
$\phi_F\colon\bC\dbl y\dbr\to
\bC\dbl x\dbr$, where $y=(y_1,\ldots, y_p)$, defined
by 
$\phi_F(f)=f\circ F$ for all $f\in\bC\dbl y\dbr$.
Conversely, any $\bC$-linear ring homomorphism
$\phi\colon\bC\dbl y\dbr\to 
\bC\dbl x\dbr$ is of the form $\phi=\phi_F$ for a
uniquely determined formal mapping 
$F\colon  (\bC^k,0)\to (\bC^p,0)$. A formal mapping
$H\colon  (\bC^k,0)\to (\bC^k,0)$ is called a {\it
formal change of coordinates} if the $k\times
k$ matrix
$\d H/\d x(0)$ is invertible. (This is equivalent to
$\phi_H$ being a
$\bC$-linear ring isomorphism.) We note that if
$\Sigma\subset\bC^k$ is a formal manifold of
codimension $d$, then there exists a formal change of
coordinates
$x'=H(x)$ so that
$\mathcal I(\Sigma)$ is generated by $x'_1,\ldots,
x'_d$. 

Let $\Lambda\subset\bC^k$, $\Sigma\subset\bC^p$ be
formal manifolds, and $F\colon (\bC^k,0)\to
(\bC^p,0)$ a formal mapping. We shall say that $F$
maps
$\Lambda$ into $\Sigma$, and write
$F(\Lambda)\subset\Sigma$, if $\phi_F(\mathcal I
(\Sigma))\subset\mathcal I(\Lambda)$.  In
particular, if
$\Lambda=\bC^k$ (so that $\mathcal I(\Lambda)=(0)$),
then
$$ F(\bC^k)\subset\Sigma \iff f\circ F=0,\ \forall
f\in\mathcal I(\Sigma).
$$ Recall that $\Rk F$ denotes the rank of the $\Bbb
K_x$-linear mapping $\Bbb K_x^k\to \Bbb K_x^p$
defined by the $p\times k$ matrix $\d F/\d x$. If
$F(\bC^k)\subset
\Sigma\subset
\bC^p$ and $\Rk F=\dim \Sigma$ (where $\dim\Sigma$
denotes the dimension of the ideal $\mathcal
I(\Sigma)$), then we shall say that $F$ is a {\it
formal subparametrization} of $\Sigma$. If, in
addition,
$\dim\Sigma=k$ and $\rk \d F/\d x(0)=k$, then we
shall say that $F$ is a {\it formal parametrization}
of
$\Sigma$. The proof of the following proposition is
elementary (using e.g.\ \cite{BER}, Proposition
5.3.5) and is left to the reader.   
\begin{Pro}\Label{fact-1} If $F\colon (\bC^k,0)\to
(\bC^p,0)$ is a formal subparametrization of a
formal manifold
$\Sigma \subset \bC^p$ then $$\mathcal I(\Sigma)=
\{f\in\bC\dbl y\dbr: f\circ F=0\}.$$
\end{Pro}

The following is a direct consequence of the formal
implicit function theorem.

\begin{Pro}\Label{fact-2} Let $F, G \colon
(\bC^k,0)\to (\bC^p,0)$ be formal mappings and
$\Sigma\subset\bC^p$ a formal manifold of dimension
$k$. If $F$ and $G$ are parametrizations of
$\Sigma$, then there exists a formal change of
coordinates
$H\colon (\bC^k,0)\to (\bC^k,0)$ such that $F\circ
H=G$.
\end{Pro}

A {\it formal vector field} $X$ in $\bC^k$ is a
$\bC$-linear derivation of the ring $\bC\dbl x\dbr$,
$x=(x_1,\ldots,x_k)$. Any formal vector field $X$ has
a unique representation in the form
$$ X=\sum_{j=1}^k a_j(x)\frac{\d}{\d x_j},\quad
a_j\in\bC\dbl x\dbr.
$$ We shall write $\D(\bC\dbl
x\dbr)$ for the $\bC\dbl
x\dbr$-module of all formal vector fields in
$\bC^k$. If
$\Sigma\subset\bC^k$ is a formal manifold, we say
that a formal vector field $X$ is {\it tangent} to
$\Sigma$ if
$Xf\in\mathcal I(\Sigma)$ for all $f\in\mathcal
I(\Sigma)$. 

For the results described in the introduction, we
shall need to focus on the situation where
$x=(Z,\zeta)$,
$Z=(Z_1,\ldots, Z_N)$ and
$\zeta=(\zeta_1,\ldots,\zeta_N)$, and consider those
formal changes of coordinates in
$\bC^N\times\bC^N$ which are of the form 
\begin{equation}\Label{H-Hbar}\mathcal
H(Z,\zeta)=(H(Z),\bar H(\zeta)),
\end{equation} where $H\colon (\bC^N,0)\to
(\bC^N,0)$ is a formal change of coordinates and
$\bar H(\zeta)$ denotes $\overline{H(\bar \zeta)}$.
We define a conjugate linear isomorphism
$\sigma\colon\bC\dbl Z,\zeta\dbr\to\bC\dbl
Z,\zeta\dbr$ by
$\sigma(f(Z,\zeta))=\bar f(\zeta, Z)$. Observe that
$\sigma^2:=\sigma\circ\sigma$ equals the identity
and that all formal changes of variables of the form
(\ref{H-Hbar}) commute with
$\sigma$, i.e.\ $\phi_{\mathcal
H}\circ\sigma=\sigma\circ\phi_{\mathcal H}$.
(Indeed, the formal changes of variables $\mathcal
H$ in
$\bC^N\times\bC^N$ of the form (\ref{H-Hbar}) are
precisely those for which $\phi_{\mathcal H}$
commute with $\sigma$ and preserve the subring
$\bC\dbl Z\dbr\subset\bC\dbl Z,\zeta\dbr$.) 

We shall say that $f\in\bC\dbl Z,\zeta\dbr$ is {\it
real} if
$\sigma(f)=f$. Hence if $f$ is a convergent power
series, then
$f$ is real if and only if $Z\mapsto f(Z,\bar Z)$ is
a real-valued function in a neighborhood of $0$ in
$\bC^N$. A manifold ideal $I\subset\bC\dbl
Z,\zeta\dbr$ is called {\it real} if
$\sigma(I)\subset I$. (Equivalently, we say that a
formal manifold $\mathcal M\subset\bC^N\times\bC^N$
is real if its ideal $\mathcal I(\scrM)$ is real.)
One can easily check that $I$ is real if and only if
there are real generators for 
$I$. The motivation for this terminology is that if
$\rho_1,\ldots,
\rho_d$ are real generators for a real manifold ideal
$I$ which are also convergent, then the equations
$\rho_1(Z,\bar Z)=\ldots\rho_d(Z,\bar Z)=0$ define a
real-analytic submanifold of codimension $d$ near
$0$ in
$\bC^N$. 

Let $\mathcal M\subset\bC^N\times\bC^N$ be a formal
real  manifold of codimension $d$ and
$\rho=(\rho_1,\ldots,\rho_d)$ a vector of generators
of $\mathcal I(M)$.  We say that $\scrM$ is {\it CR}
if
$$\rk_{\Bbb K_{Z,\zeta}}\frac {\d \rho}{\d Z}
(Z,\zeta)=
\rk
\frac {\d \rho}{\d Z}
(0,0).$$
Here $\d \rho/\d Z$ is considered as a $d\times N$
matrix with coefficients in $\bC\dbl Z,\zeta\dbr$ and
$\Bbb K_{Z,\zeta}$ is the field of fractions of
that ring.
We say that the formal real submanifold 
$\scrM\subset \bC^N$ is {\it generic} if $\rk
\d \rho/\d Z(0)=d$. Observe that the definitions
above are independent of the choice of generators
$\rho_j$, and that any formal generic manifold is
necessarily CR.

 If
$\scrM\subset\bC^N\times\bC^N$ is a formal CR
manifold, then it can be shown (cf.\ Remark
\ref{CR-ideal}) that
$\mathcal I(\scrM)\cap\bC\dbl Z\dbr$ is a manifold
ideal in
$\bC\dbl Z\dbr$ of dimension 
$N-d+\rk \d\rho/\d Z(0)$. The formal manifold
$\mathcal C\subset\bC^N$ whose ideal is 
$\mathcal I(\scrM)\cap\bC\dbl Z\dbr$ is called the
{\it intrinsic complexification} of
$\mathcal M$. Observe that $\scrM$ is generic if and
only if
$\mathcal C=\bC^N$.

A {\it formal $(1,0)$-vector field} is a formal
vector field in $\bC^N\times\bC^N$ of the form 
$$ X=\sum_{j+1}^N a_j(Z,\zeta)\frac{\d}{\d Z_j},\quad
a_j\in\bC\dbl Z,\zeta\dbr.
$$  Similarly, a {\it formal $(0,1)$-vector field}
is one of the form
$$ X=\sum_{j+1}^N a_j(Z,\zeta)\frac{\d}{\d
\zeta_j},\quad a_j\in\bC\dbl Z,\zeta\dbr.
$$
If $\scrM\subset\bC^N\times\bC^N$ is a formal CR
manifold, then a formal $(0,1)$-vector
field tangent to $\scrM$ is called a {\it CR vector
field } on $\scrM$. We denote by $\frak g_{\scrM}$
the Lie algebra of formal vector fields generated by
all the formal $(0,1)$ and $(1,0)$-vector fields
tangent to $\scrM$. One can easily check that $\frak
g_{\scrM}$ is also a $\bC\dbl Z,\zeta\dbr$-module. If
$\dim\frak g_{\scrM}(0)=\dim\scrM$, then we say that
$\scrM$ is of {\it finite type}. (Here $\frak
g_{\scrM}(0)$ is the complex vector space spanned by
the vector fields in $\frak
g_{\scrM}$ evaluated at $0$.)

We extend the domain of definition of the involution
$\sigma$ to include formal vector fields as follows
$$
\sigma(X)f:=\sigma(X(\sigma(f))),\quad f\in\bC\dbl
Z,\zeta\dbr.
$$ It is easy to check that if
$X=\sum_{\j=1}^Na_j\d/\d Z_j+b_j\d/
\d\zeta_j$, then
$$\sigma(X)=\sum_{\j=1}^N\sigma(a_j)\d/\d
\zeta_j+\sigma(b_j)\d/
\d Z_j.$$ In particular, if $X$ is a formal
$(0,1)$-vector field, then
$\sigma(X)$ is a formal $(1,0)$-vector field (and
conversely). Moreover, if
$\scrM\subset\bC^N\times\bC^N$ is a formal real
manifold and
$X$ is a formal vector field tangent to $\scrM$, then
$\sigma(X)$ is tangent to $\scrM$. Hence, if $\scrM$
is a formal CR manifold, and $\mathcal
O\subset\scrM$ is a formal real manifold to which
all the formal CR vector fields on $\scrM$ are
tangent, then all vector fields in
$\frak g_{\scrM}$ are also tangent to $\mathcal O$.

If $\scrM\subset\bC^N\times\bC^N$ is a formal
generic manifold, then we define the Segre variety
mapping
$\gamma$ and the iterated Segre mappings $v^j$
exactly as given by (\ref{gamma}) and (\ref{vj}) in
the introduction.  If $M\subset \bC^N$ is a smooth
generic submanifold through the origin, then we
associate to it a formal generic manifold
$\scrM\subset \bC^N\times\bC^N$ as follows. Let
$r_1,\ldots, r_d$ be smooth defining functions for
$M$ near
$0$ and $\rho_1(Z,\bar Z),\ldots,\rho_d(Z,\bar Z)$
their Taylor series at $0$. The formal generic
manifold $\scrM\subset\bC^N\times\bC^N$ is defined
to be that associated to the ideal $\mathcal
I(\scrM)\subset \bC\dbl Z,\zeta\dbr$ generated by
the formal power series
$\rho_1(Z,\zeta),\ldots,\rho_d(Z,\zeta)$. Observe
that $\frak g_{\scrM}$ consists of the
formal vector fields obtained taking the Taylor
expansions of the coefficients of the smooth vector
fields in
$\frak g_M$. In particular, $\dim \frak
g_{\scrM}(0)=\dim \frak g_M(0)$. Hence, Theorems
\ref{main}, \ref{main-2}, and \ref{main-3} are
consequences of the following more general result.

 \begin{Thm}\Label{formal-main} Let $\scrM\subset
\bC^N\times\bC^N$ be a formal generic manifold of
codimension $d$, and let
$\gamma(\zeta,t)$ be a Segre variety mapping of
$\scrM$. Let $v^j$, $j\geq 1$, be the iterated Segre
mappings of
$\scrM$ relative to $\gamma$. Then, the rank $\Rk
v^j$ is an increasing function of $j$ and is
independent of the choice of the formal holomorphic
coordinates
$Z$, the formal power series $\rho_j$, and the Segre
variety mapping $\gamma$.  In addition, there exists
an integer $k_0$,
$1\le k_0\leq d+1$, such that $\Rk v^j=\Rk v^{j+1}$
for
$j\ge k_0$, and if $k_0> 1$, then $\Rk v^j < \Rk
v^{j+1}$ for $1\le j\le k_0-1$. Moreover, the
following holds.
\begin{enumerate}
\item [(i)] There exists a unique formal CR
manifold
$\mathcal O\subset \scrM$, with $\dim_\bR \mathcal
O=\dim_\bC\frak g_{\scrM}(0)$, such that all the
formal CR vector fields on
$\scrM$ are tangent to $\mathcal O$. Moreover, if
$e=2N-d-\dim\frak g_{\scrM}(0)$, then there exist
$f_1,\ldots, f_e\in\bC\dbl Z\dbr $ such that
$f_1,\ldots, f_e,\rho_1,\ldots,\rho_d$ generate the
manifold ideal
$\mathcal I(\mathcal O)$. 
\item[(ii)] Let $\mathcal W\subset\bC^N$ be the
formal complex manifold defined by the ideal
$I(f_1,\ldots f_e)\subset\bC\dbl Z\dbr$. Then 
$v^j(\bC^{jn})\subset \mathcal W$ for
all $j\ge 1$.
\item [(iii)] $\Rk v^{k_0}=\dim\frak
g_{\scrM}(0)+d-N=\dim\mathcal W$.
\end{enumerate}
\end{Thm}

We should point out that, in general, the rank of
the $N \times jn$ matrix
$$ 
 \frac{\d
v^{j}}{\d(t^1,\ldots,t^{j})}(0)
$$ is strictly less than $\Rk
v^{k_0}$ for all $j$. However, we have the
following.

\begin{Thm}\Label{formal-double} Let
$\scrM\subset\bC^N\times\bC^N$ be a formal generic
manifold of codimension $d$, and let $\gamma$,
$v^j$, and $k_0$ be as in Theorem
$\ref{formal-main}$. Then there exists a formal
manifold $\Sigma\subset\bC^{2nk_0}$, with $n=N-d$,
of dimension $nk_0$ such that
$$v^{2k_0}(\Sigma)=\{0\},$$ and for any formal
parametrization
$F\colon (\bC^{nk_0},0)\to (\bC^{2nk_0},0)$ of
$\Sigma$,
$$ 
\rk_{\Bbb K_s}\left( \frac{\d
v^{2k_0}}{\d(t^1,\ldots,t^{2k_0})}(
F(s^1,\ldots,s^{ k_0}))\right)=\Rk v^{k_0},
$$ where $s=(s^1,\ldots,s^{ k_0})\in \bC^{nk_0}$.
\end{Thm}

In the case where $M\subset\bC^N$ is a real-analytic
submanifold, say of finite type at $0$, Theorem
\ref{formal-main} asserts that the image of
$v^{k_0}$ contains an open subset of $\bC^N$, but not
necessarily an open neighborhood of $0$. However,
Theorem
\ref{formal-double} asserts that the image of
$v^{2k_0}$ contains an open neighborhood of $0$ in
$\bC^N$. 

The proof of Theorem \ref{formal-main} will be given
in \S\S \ref{pfmain-1}--\ref{pfmain-6}. The proof
of Theorem
\ref{formal-double} is given in \S
\ref{pfmain-6}.

\section{Increase of the rank of the iterated Segre
mappings}\Label{pfmain-1}
In this section we begin the proof of Theorem
\ref{formal-main}.  We
start by showing that
$\Rk v^j$ is an increasing function of $j$, i.e.\
$\Rk v^{j+1}\geq
\Rk v^j$. We shall restrict ourselves to the case
where
$j$ is odd; the even case is similar and left to the
reader. By iterating the definition (\ref{vj}) of
$v^j$, we obtain
\begin{equation}\Label{vj-it}
v^{j+1}(t^1,\ldots,t^{j+1})=\gamma(\bar\gamma(\ldots
\gamma(\bar\gamma(0,t^1),t^2),\ldots,t^j),t^{j+1}).
\end{equation} Hence, setting $t^1=0$ and using the
fact that
$\gamma(0,0)=0$, we
have
\begin{equation}\Label{4}
\begin{aligned} v^{j+1}(0,t^2,\ldots,t^{j+1})= &\,
\gamma(\bar\gamma (\ldots
\gamma(\bar\gamma(0,0),t^2),\ldots,t^j),t^{j+1})\\
=&\, v^j(t^2,\ldots,t^{j+1}).
\end{aligned}
\end{equation} The conclusion $\Rk v^{j+1}\geq \Rk
v^j$ follows immediately from (\ref{4}).

We shall now prove that if, for some integer $k_0$,
$\Rk  v^{k^0}=\Rk v^{k_0+1}$, then $\Rk v^{j}=\Rk
v^{k_0}$ for all
$j\geq k_0$. In view of (\ref{vj}) and  (\ref{4}),
this statement is a consequence of the following
proposition.

\begin{Pro}\Label{prop-er} Let
$A\colon(\bC^r\times\bC^s,0)\to (\bC^m,0)$ and
$F\colon (\bC^n\times\bC^m,0)\to (\bC^N,0)$ be
formal mappings, and set $B(x)=A(x,0)$, where
$x=(x_1,\ldots, x_r)$. If
$\Rk A(x,y)=\Rk B(x)$, then $\Rk F(z,A(x,y))=\Rk
F(z,B(x))$. Here, $z=(z_1,\ldots, z_n)$ and
$y=(y_1,\ldots, y_s)$.
\end{Pro}

\begin{proof} Let $B^j(x,y_1,\ldots,
y_j)=A(x,y_1,\ldots, y_j,0,\ldots, 0)$, for
$j=1,\ldots, s-1$. Clearly, if $\Rk A(x,y)=\Rk
B(x)$, then $\Rk A(x,y)=\Rk B^j(x,y_1,\ldots, y_j)$
for each $j=1,\ldots, s-1$. A simple induction
argument shows that it suffices to prove the
proposition when $s=1$, which we shall assume for
the remainder of the proof of Proposition
\ref{prop-er}. The Proposition will be a consequence
of  the following lemma. 
\begin{Lem}\Label{lemma-flow}Let $A$, $B$ be as in
Proposition
$\ref{prop-er}$ with $s=1$. If $\Rk A(x,y)=\Rk
B(x)$, then there exists a formal mapping
$\phi\colon (\bC^r\times\bC,0)\to
(\bC^r\times\bC,0)$ such that $\Rk \phi(u,t)=r+1$,
where
$u=(u_1,\ldots, u_r)$, and 
\begin{equation}\Label{flow-1}
A(\phi(u,t))=A(\phi(u,0))=B(u).
\end{equation}
\end{Lem}
We now show that Proposition
\ref{prop-er} follows from Lemma \ref{lemma-flow},
whose proof will be given below. Set
$H(z,x,y)=F(z,A(x,y))$ and 
$G(z,u,t)=F(z,A(\phi(u,t)))$. Observe that the
formal mapping
$\psi\colon (\bC^n\times\bC^r\times\bC,0)\to
(\bC^n\times\bC^r\times\bC,0)$, where
$\psi(z,u,t)=(z,\phi(u,t))$ satisfies $\Rk
\psi(z,u,t)=n+r+1$. Hence, since $G=H\circ \psi$, it
follows that $\Rk H(z,x,y)=\Rk G(z,u,t)$. Also, by
Lemma \ref{lemma-flow}, $G(z,u,t)=
F(z,A(\phi(u,t)))=F(z,B(u))$, which proves
Proposition
\ref{prop-er}. 
\end{proof}

\begin{proof}[Proof of Lemma $\ref{lemma-flow}$] Let
$p=\Rk A(x,y)=\Rk B(x)$, and hence $p \le r$.  
Without loss of generality, we may assume that
$\partial B/\partial x_j\in \bC\dbl x\dbr^m$,
$j=1,\ldots p$, are linearly independent over $\Bbb
K_x$ and hence $\partial A/\partial x_j\in \bC\dbl
x,y\dbr^m$,
$j=1,\ldots p$, are also linearly independent over
$\Bbb K_{x,y}$. We write $x=(x',x'')$, with
$x'=(x_1,\ldots, x_p)$, and e.g.
$\partial A/\partial x'$ for the $m\times p$
Jacobian matrix of $A$ with respect to $x'$.
Consider the $\Bbb
K_{x,y}$-linear system
\begin{equation}\Label{system}
\pmatrix \partial A/\partial x' & \partial
A/\partial x''& \partial A/\partial y\endpmatrix
\pmatrix a' \\ a''\\b\endpmatrix =\pmatrix
0\\0\\0\endpmatrix
\end{equation}
with $a'\in \Bbb K_{x,y}^p$, $a''\in \Bbb
K_{x,y}^{r-p}$, and
$b\in \Bbb K_{x,y}$.  We take
$a''(x,y)=0$, $b(x,y)=1$ and solve (\ref{system})
for 
$a'$, which is possible since the rank of $\partial
A/\partial x'$ is
$p$. Moreover, if we denote by $k_j$, $j=1,\ldots,
p$, indices such that $\Delta:=\det(\partial
A_{k_j}/\partial x_i)_{1\leq i,j\leq p}$ with
$\Delta (x,0)\neq 0$, then the solution $a'$ is of the
form
\begin{equation}\Label{c}
a'(x,y)=\frac{c(x,y)}{\Delta(x,y)},
\end{equation} where $c\in\bC\dbl x,y\dbr^p$. Write
\begin{equation}\Label {vectorf}
\mathcal X=\Delta(x,y)\frac{\d}{\d y}+\sum_{j=1}^p
c_j(x,y)\frac{\d}{\d x_j},
\end{equation} where $c=(c_1,\ldots, c_p)$ is given
by (\ref{c}). Observe that the coefficients of
$\mathcal X$ are in $\bC\dbl x,y\dbr$, 
$\mathcal X A(x,y)=0$, and $\Delta(x,0)\neq 0$. We
define
$\phi(u,t)$, $u = (u_1, \ldots, u_r)$, to be the
formal flow of
$\mathcal X$, i.e.\ the solution of the initial value
problem
$$
\frac{\d\phi}{\d t}(u,t)=X(\phi(u,t)),\quad
\phi(u,0)=(u,0),
$$ where $X(x,y)$ denotes the vector of coefficients
of $\mathcal X$. The fact that $\Delta(x,0)\neq 0$
implies that $\Rk
\phi(u,t)=r+1$. Since $\mathcal X A=0$, it is easy
to verify that the identity (\ref{flow-1}) holds.
This completes the proof of Lemma \ref{lemma-flow}
and hence that of Proposition \ref{prop-er}.
\end{proof}

In summary, we have shown in this section that
\begin{equation}\Label{incr}
n = \Rk v^1 \le \Rk v^j \le \Rk v^{j+1} \le N, \ \ 
\forall j \ge 1.
\end{equation}
Moreover, there exists an integer $k_0$, $1 \le k_0
\le d+1$, such that
\begin{equation}\Label{strincr}
\Rk v^j < \Rk v^{j+1},  \  1 \le j < k_0, \ \ {\rm
and}\ \ \Rk v^{k_0} = \Rk v^j, \ \ j \ge k_0.
\end{equation}

\section{Independence of rank on choice of Segre
variety mapping}\Label{pfmain-2}

In view of (\ref{incr}) and (\ref{strincr}), to
complete  the proof of the first part of Theorem
\ref{formal-main}, it suffices to show that
$\Rk v^j$ is independent of the choice of the Segre
variety mapping
$\gamma$. For this, we introduce, for each positive
even integer $2m$, the formal manifold $\frak
S_{2m}\subset\bC^{2mN}$ whose ideal in $\bC\dbl
Z,\zeta^1,Z^1,\ldots,Z^{m-1},\zeta^m\dbr$, where
$Z=(Z_1,\ldots,Z_N)$,
$Z^j=(Z^j_1,\ldots Z^j_N)$,
$\zeta^j=(\zeta^j_1,\ldots,\zeta^j_N)$, is generated
by  
\begin{equation} \Label{ind-1}
\rho(Z,\zeta^1),\ \rho(Z^1,\zeta^1),\
\rho(Z^1,\zeta^2),\
\rho(Z^2,\zeta^2),\ \ldots,\rho(Z^{m-1},\zeta^{m}),\
\rho(0,\zeta^{m}),
\end{equation} where $\rho=(\rho_1,\ldots,\rho_d)$
are formal defining functions of
$\scrM$.  Similarly, for positive odd integers
$2m+1$, we define $\frak
S_{2m+1}\subset\bC^{(2m+1)N}$ to be the formal
manifold  whose ideal in $\bC\dbl
Z,\zeta^1,Z^1,\ldots,Z^{m-1},\zeta^m,Z^m\dbr$ is
generated by  
\begin{equation} \Label{ind-2}
\rho(Z,\zeta^1),\ \rho(Z^1,\zeta^1),\
\rho(Z^1,\zeta^2),\
\ldots,\rho(Z^{m},\zeta^{m}),\ \rho(Z^{m},0). 
\end{equation} The formal manifold $\frak
S_k\subset\bC^{kN}$ has dimension $kn$. It is called
the $k$th  Segre manifold of $\scrM$. Also,
denote by
$\pi^k_1\colon(\bC^{kN},0)\to (\bC^N,0)$ the  
projection on the first factor, e.g.
$$
\pi^{2m+1}_1(Z,\zeta^1,Z^1,\ldots,\zeta^m,Z^m)=Z.
$$ Fix a choice of Segre variety mapping $\gamma$. 
Define, for positive even integers $2m$, the formal
mapping
$T^{2m}\colon (\bC^{2mn},0)\to (\bC^{2mN},0)$ by
\begin{multline}\Label{ind-3}
T^{2m}(t^1,\ldots,t^{2m})\colon=\\\big(v^{2m}(t^1,\ldots,
t^{2m}),\bar v^{2m-1}(t^1,\ldots,
t^{2m-1}),\ldots,\bar v^1(t^1)\big),\
\end{multline}  and
$t^j=(t^j_1,\ldots, t^j_{n})$.  Similarly, for
positive odd integers $2m+1$, define the formal
mapping
$T^{2m+1}\colon (\bC^{(2m+1)n},0)\to
(\bC^{(2m+1)N},0)$ by 
\begin{multline}\Label{ind-4}
T^{2m+1}(t^1,\ldots,t^{2m+1})\colon
=\\\big(v^{2m+1}(t^1,\ldots, t^{2m+1}),\bar
v^{2m}(t^1,\ldots, t^{2m}),
\ldots,v^1(t^1)\big),\
\end{multline} Observe, by (\ref{vj}) and the
definition (\ref{gamma}) of the Segre variety
mapping, that each $T^k$ in fact maps $\bC^{kn}$
into $\frak S_k$, i.e. $f\circ T^{k}= 0$ for each
$f$ in the ideal of
$\frak S_k$. Moreover, as is easy to verify from the
fact that the rank of
$\d \gamma/\d t(0)$ is $n$, the rank of $\d T^k/\d
(t^1,\ldots, t^k)(0)$ equals $kn$ which is also the
dimension of $\frak S_k$. Hence,
$T^k$ is a parametrization of the formal manifold
$\frak S_k$.  In addition, we have 
$v^k=\pi^k_1\circ T^k$. 

Now, denote by
$\tilde \gamma$ another choice of Segre variety
mapping, and by
$\tilde v^k$, $\tilde T^k$ the corresponding
mappings defined as above using $\tilde\gamma$
instead of $\gamma$. For the same reasons as above,
the mapping $\tilde T^k\colon(\bC^{kn},0)\to
(\bC^{kN},0)$ is a parametrization of the formal
manifold $\frak S_k$, and $\tilde
v^k=\pi^k_1\circ\tilde T^k$. Since both
$T^k\colon (\bC^{kn},0)\to (\bC^{kN},0)$ and $\tilde
T^k\colon (\bC^{kn},0)\to (\bC^{kN},0)$  are
parametrizations of $\frak  S_k$, there exists a
formal invertible mapping
$F^k\colon\:(\bC^{kn},0)\to (\bC^{kn},0)$ such that
$\tilde T^k=T^k\circ F^k$. Thus,
$\tilde v^k=v^k\circ F^k$ and the fact that $\Rk
\tilde v^k=\Rk v^k$ is a consequence of the chain
rule. This completes the  proof of the statement in
Theorem \ref{formal-main} that $\Rk v^k$ is
independent of the choice of
$\gamma$. 

\section{Construction and properties of the formal CR
manifold $\mathcal O$}\Label{pfmain-3}

To construct the formal CR manifold $\mathcal
O\subset\mathcal M$ in Theorem
\ref{formal-main} (i), we shall need the following,
which may be regarded as a formal version of
Nagano's theorem (see e.g. \cite {BER}, Theorem
3.1.4) for real analytic vector fields. As mentioned 
in the introduction, the inductive proof given here is based
on a suggestion by D. Zaitsev.

\begin{Pro}\Label{Nagano}Let $\frak g$  be a Lie
algebra of formal vector fields in $k$
indeterminates $x = (x_1,\ldots, x_k)$ and assume
that $\frak g$ is also a
$\bC[[x]]$ module.  Then there exists a unique formal
manifold $\Sigma \subset \bC^k$ with
$\dim_{\bC}\Sigma=\dim_{\bC}\frak g(0)$ such that
any $X \in
\frak g$ is tangent to $\Sigma$. Moreover, if
$\Lambda\subset\bC^k$ is any other formal manifold
such that all formal vector fields in $\frak g$ are
tangent to
$\Lambda$, then $\Sigma\subset\Lambda$.
\end{Pro}

\begin{proof} We shall prove the proposition by
induction on $k$.  If
$k=1$, then the proposition holds with either
$\Sigma =
\{0\}$ or $\Sigma = \bC$, depending on whether or
not all vector fields in $\frak g$ vanish at $0$. 
Now assume that the proposition holds for $k$;  we
shall prove it for $k+1$.  If all the $X \in \frak
g$ vanish at
$0$, then again $\Sigma = \{0\}$ is the only formal
manifold satisfying the conclusion of the
proposition. If not, without loss of generality, and
after making a formal change of coordinates, we may
assume that  the vector field
$X_1 = {\partial\over \partial x_1}$ is in $\frak g$.

Consider the Lie subalgebra $\frak g'\subset \frak g$
of all formal vector fields in $\frak g$ of the form
$\sum_{j=2}^{k+1} a_j(x){\partial\over\partial
x_j}$.  It is easy to see that any $X \in \frak g$
may be written uniquely in the form $X = a(x)X_1 +
X'$ with $a(x) \in
\bC\dbl x \dbr$ and $X'\in \frak g'$.  We write $x'
= (x_2,\ldots, x_{k+1})$ and consider the Lie
algebra $\frak g''$ of formal vector fields in the
indeterminates $x'$ obtained from $\frak g'$ by
replacing
$x_1$ by $0$. By the inductive hypothesis, there is a
formal manifold $\Sigma'\subset \bC^k$ to which
all the vector fields in $\frak g''$ are tangent and
such that
$\dim \Sigma' = \dim_{\bC}\frak g''(0)$. Let
$w_i(x') \in
\bC\dbl x' \dbr$, $i = 1,\ldots, r$, be generators of
the manifold ideal $\mathcal I(\Sigma')$.

\begin{Lem}\Label{Nagano-lem}  Let $\Sigma\subset
\bC^{k+1}$ be the formal manifold whose manifold
ideal is generated by the 
$w_i(x')$, $i = 1,\ldots, r$, regarded as elements of
$\bC\dbl x \dbr$  (independent of $x_1$). Then all
the vector fields in $\frak g$ are tangent to
$\Sigma$ and $ r=k+1-\dim \frak g(0)$. 
\end{Lem} 
\begin{proof} [Proof of Lemma $\ref{Nagano-lem}$] We
shall first show that any
$X
\in
\frak g$ is tangent to $\Sigma$.  For $X\in \frak g$
and $X_1 = {\partial\over \partial x_1}$ as before,
we shall use the notation
$(\ad X_1) X = [X_1, X]$. Since any $X \in \frak g$
can be written in the form $X = a(x)X_1 + X'$ with
$X'\in \frak g'$, and since the $w_j$ are
independent of $x_1$, it suffices to show that any
vector field in $\frak g'$ is tangent to $\Sigma$. 
For this we expand
$(X'w_j)(x_1,x')$ as a power series in $x_1$ and
obtain
\begin{equation}\Label{ad-1}
(X'w_j)(x_1,x')=\sum_{p=0}^\infty (Y_p
w_j)(x')\frac{x_1^p}{p!},
\end{equation} where $Y_p$ is the vector field in
$\frak g''$ given by
$$ Y_p:=((\ad X_1)^p X')|_{x_1=0}.
$$ Since $Y_p\in\frak g''$ it is tangent to
$\Sigma'$ and, hence, we have 
$$ (Y_pw_j)(x')=\sum_{i=1}^r c_{ijp}(x')w_i(x'), 
$$ for some $c_{ijp}\in\bC\dbl x''\dbr$.  Hence, $X'$
is tangent to $\Sigma$ by (\ref{ad-1}). It remains
to show that
$r=k+1-\dim_\bC\frak g(0)$. Since all the vector
fields in
$\frak g$ are tangent to $\Sigma$, we have $r\leq
k+1-\dim_\bC\frak g(0)$. The opposite inequality
follows from the assumption that $r=k-\dim_\bC\frak
g''(0)$ and the fact that $\dim_\bC\frak g''(0)\leq
\dim_\bC\frak g(0)-1$. 
\end{proof}

To finish the proof of Proposition \ref{Nagano}, we
must show that the formal manifold $\Sigma$ provided
by Lemma
\ref{Nagano-lem} is unique, and if
$\Lambda\subset\bC^k$ is as in the statement of
Proposition \ref{Nagano} then
$\Sigma\subset\Lambda$.  Let $f_1,\ldots,
f_m\in\bC\dbl x_1,x'\dbr$ be generators of the
manifold ideal $\mathcal I(\Lambda)$. Since $X_1$ is
tangent to $\Lambda$, it follows that
$X_1^jf_l\in\mathcal I(\Lambda)$ for $l=1,\ldots, m$
and $j=1,2,\ldots$. Let $I''$ be the ideal in
$\bC\dbl x'\dbr$ generated by $f_1(0,x'),\ldots
f_m(0,x')$. Note that
$X_1f_l(0,0)=\d f_l/\d x_1(0,0)=0$, $l=1,\ldots, m$,
and hence
$I''$ is a manifold ideal of the same codimension as
$\mathcal I(\Lambda)$. Since
$\frak g'\subset\frak g$ is tangent to
$\mathcal I(\Lambda)$,
$\frak g''$ is tangent to $I''$. By the induction
hypothesis, the manifold ideal $\mathcal I(\Sigma')$
is unique and
$I''\subset\mathcal I(\Sigma')$. Thus, $X_1^j
f_l\in\mathcal I(\Lambda)$ implies that there are
$c_{ijl}\in\bC\dbl x'\dbr$ such that
\begin{equation}\Label{nag-1} X_1^j
f_l(0,x')=\sum_{i=1}^r c_{ijl}(x') w_i(x').
\end{equation} Now, since 
\begin{equation}\Label{nag-2}
f_l(x_1,x')=\sum_{j=0}^\infty X_1^j
f_l(0,x')\frac{x_1^j}{j!}
\end{equation} we conclude, by substituting
(\ref{nag-1}) in (\ref{nag-2}), that $f_l\in\mathcal
I(\Sigma)$ for
$l=1,\ldots, m$. Hence, $\mathcal
I(\Lambda)\subset\mathcal I(\Sigma)$. This also
proves the uniqueness of $\Sigma$, since if the
dimensions of $\Lambda$ and
$\Sigma$ were the same, then necessarily
$\Lambda=\Sigma$. This completes the proof of
Proposition \ref{Nagano}. 
\end{proof}

Let $\mathcal O\subset \bC^{N}\times\bC^{N}$ be the
formal manifold obtained by applying Proposition
\ref{Nagano} to the Lie algebra (and $\bC\dbl
Z,\zeta\dbr$ module)
$\frak g_{\scrM}$. Observe that $\mathcal O\subset
\mathcal M$ since all the vector fields in $\frak
g_M$ are tangent to
$\mathcal M$. We shall call $\mathcal O$ the {\it
formal CR orbit} of $0$ in $\mathcal M$. (We shall
show in  \S\ref{pfmain-6} that in fact
$\mathcal O$ is a formal real CR manifold.) 

We shall continue the proof of Theorem
\ref{formal-main}. Since
$\partial_Z\rho_1(0,0),\ldots\d_Z\rho_d(0,0)$ are
linearly independent, we may assume after
renumbering the coordinates, that
$Z=(z,w)$, with $z=(z_1,\ldots,z_n)$,
$w=(w_1,\ldots, w_d)$, and the $d\times d$ matrix
$\d\rho/\d w (0,0)$ is invertible. Hence, by the
formal implicit function theorem, the ideal of
$\mathcal M$  is generated by 
$w_j-Q_j(z,\Z)$,
$j=1,\ldots, d$, where the
$Q_j(z,\Z)$ are formal power series in $n+N$
indeterminates without constant terms. If we write
$\zeta=(\chi,\tau)$ with
$\chi=(\chi_1,\ldots,\chi_n)$,
$\tau=(\tau_1,\ldots,\tau_d)$,  then the reality of
$\mathcal M$ implies that we may also take
$\tau_j-\bar Q_j(\chi,z,w)$,
$j=1,\ldots, d$, as generators of the ideal $\mathcal
I(\mathcal M)$. Consequently we also have the
identity
\begin{equation}\Label{reality} Q\big(z,\chi,\bar
Q(\chi,z,w)\big)= w,
\end{equation} where $Q=(Q_1,\ldots, Q_d)$.

  As
a basis for the
$(0,1)$ and $(1,0)$ vector fields tangent to
$\mathcal M$, we choose
\begin{equation}\Label{Lj}
\begin{aligned}
\mathcal L_j:= &\,\frac{\d}{\d \chi_j}+\sum
_{l=1}^d\bar Q_{l,\chi_j}(\chi,z,w)\frac{\d}{\d
\tau_l},\quad j=1,\ldots,n,\\ \tilde \mathcal L_j:=
&\,\frac{\d}{\d z_j}+\sum _{l=1}^d
Q_{l,z_j}(z,\chi,\tau)\frac{\d}{\d w_l},\quad
j=1,\ldots,n.
\end{aligned}
\end{equation}

\begin{Pro}\Label{orbit-ideal} Let $\mathcal
M\subset\bC^N$ be a formal generic manifold whose
ideal is generated by
$\rho_1,\ldots,\rho_d\in\bC\dbl Z,\zeta\dbr$. Let
$\mathcal O$ be the formal CR orbit of $0$ in
$\mathcal M$ and
$e:=2N-d-\dim_\bC\mathcal O$. Then there are
$f_1,\ldots,f_e\in\bC\dbl Z\dbr$ such that the formal
manifold ideal $\mathcal I(\mathcal O)$ is generated
by
$\rho_1,\ldots, \rho_d,f_1,\ldots, f_e$.
\end{Pro}

\begin{proof} It suffices to prove the proposition
using the special choice of coordinates $Z=(z,w)$,
$\zeta=(\chi,\tau)$ made above. Thus, we may assume
that
$\rho_j(Z,\zeta)=\tau_j-\bar Q_j(\chi,z,w)$,
$j=1,\ldots, d$. Since $\mathcal O$ is a formal
manifold of codimension
$d+m$ and $\mathcal I(\mathcal M)\subset\mathcal
I(\mathcal O)$, there are
$g_1,\ldots, g_e\in\bC\dbl Z,\zeta\dbr$ such that
$\mathcal I(\mathcal O)$ is generated by
$\rho_1,\ldots, \rho_d,g_1,\ldots, g_e$. By using the
special form of the $\rho_j$ (i.e.\ substituting
$\bar Q_j(\chi,z,w)$ for $\tau_j$), we may assume
that
$g_j$ is independent of $\tau$, that is,
$g_j=g_j(z,w,\chi)$, $j=1,\ldots, e$. Since the
$\mathcal L_j$ are tangent to $\mathcal O$ and
$\mathcal L_j g_l=\d g_l/\d \chi_j$, $1\leq l\leq e$
and $1\leq j\leq n$, we obtain 
\begin{equation}\Label{de-1}
\frac{\d
g}{\d\chi_l}(z,w,\chi)=A_l(z,w,\chi)g(z,w,\chi),\quad
l=1,\ldots, n,
\end{equation} where $g=(g_1,\ldots, g_e)$ and the
$A_l$ are $(e\times e)$-matrices whose entries are
in $\bC\dbl z,w,\chi\dbr$. Let $u^1,\ldots,
u^e\in\bC\dbl z,w,\chi\dbr^e$ be a fundamental
system of solutions for the system of differential
equations given by (\ref{de-1}) with
$l=1$. Denote by $U(z,w,\chi)$ the $(e\times
e)$-matrix in which the $u^j$ are columns. Since $g$
is a solution of the system (\ref{de-1}) (in
particular with
$l=1$) there exists
$c^1\in\bC\dbl z,w,\chi'\dbr^e$, where
$\chi=(\chi_1,\chi')$, such that $g=Uc^1$. Since $U$
is invertible, we have
$c^1(z,w,\chi')=(U(z,w,\chi))^{-1}g(z,w,\chi)$ and,
hence, each component of $c^1(z,w,\chi')$ is in
$\mathcal I(\mathcal O)$. Moreover, $\rho_1,\ldots,
\rho_d, c^1_1,\ldots, c^1_e$ generate $\mathcal
I(\mathcal {O})$. Proceeding this way to eliminate
$\chi_2,\ldots, \chi_n$, we finally obtain
$c^n\in\bC\dbl z,w\dbr^e$ such that
$f_l(z,w)=c^n_l(z,w)$, $l = 1,\ldots, e$, satisfy the
conclusion of Proposition
\ref{orbit-ideal}.  
\end{proof}

\begin{Rem}\Label{CR-ideal} {\rm By combining the
linear algebra argument in the proof of \cite{BER},
Theorem 1.8.1, with the argument in the proof of
Proposition \ref{orbit-ideal} above, one can show
the following. If $\rho_1,\ldots, \rho_d$ generate
the ideal $\mathcal I(\scrM)$ of a formal CR manifold
$\scrM\subset\bC^N\times\bC^N$, then there are
$\tilde\rho_1,\ldots,\tilde\rho_{d-e}\in\bC\dbl
Z,\zeta\dbr$, $f_1,\ldots, f_e\in\bC\dbl Z\dbr$,
where
$e=d-\rk\d\rho/\d Z(0)$, such that
$\tilde\rho_1,\ldots,\tilde\rho_{d-e}, f_1,\ldots,
f_e$ also generate $\mathcal I(\scrM)$. Moreover,
the rank of $\d
\tilde\rho/\d Z(0)$ equals the rank of $\d
\rho/\d Z(0)$ and, hence, the rank of $\d f/\d Z (0)$
is
$e$. The formal manifold $\mathcal C\subset\bC^N$
whose ideal is generated by $f_1,\ldots, f_e$ is the
intrinsic complexification of $\mathcal M$ defined
in \S
\ref{Prelim}. 
 }
\end{Rem}

We have now proved the properties of $\mathcal O$
announced in Theorem \ref{formal-main} (i) with the
exception of the claim that $\mathcal O$ is CR. This
will be done in \S
\ref{pfmain-6}.

\section{Rank of mappings
$\Theta^j$, $\Phi^j$ into $\mathcal
O$}\Label{pfmain-4}

To prove statement (iii) of Theorem
\ref{formal-main}, we shall construct special
mappings $\Theta^j$, $\Phi^j$ into
$\mathcal O$ and use these mappings to compute the
rank of
$v^j$. We shall proceed using a
special choice of Segre variety mapping $\gamma$.
We use the formal coordinates $Z=(z,w)$ and $\zeta =
(\chi, \tau)$ introduced in
\S
\ref {pfmain-3}. Given any Segre variety mapping
$\gamma(\Z , t)$, we may decompose it as
$$
\gamma(\zeta,t)=(\mu(\Z,t),\nu(\Z,t)),
$$ with $\mu=(\mu_1,\ldots,\mu_n)$ and
$\nu=(\nu_1,\ldots,\nu_d)$. It follows that 
$$
\nu(\Z,t)=Q(\mu(\Z,t),\Z),
$$ and hence necessarily $\rk (\d \mu/\d t (0,0))=n$.
We now make the choice of $\gamma(\zeta,t)$
corresponding to $\mu(\zeta,t)=t$, i.e.
\begin{equation}\Label{gamma-spec}
\gamma(\Z,t)=(t,Q(t,\zeta))=(t, Q(t,\chi,\tau)).
\end{equation} We need the following lemma.

\begin{Lem}\Label{vj+1=vj-1} Let $\gamma$ be given by
  $(\ref{gamma-spec})$. Then for any integer $j\geq
2$, the following holds:
\begin{equation}\Label{hj}
v^{j+1}(t^1,\ldots,t^{j-1},t^j,t^{j+1})\big|_{t^{j+1}
=t^{j-1}} = v^{j-1}(t^1,\ldots,t^{j-1}).
\end{equation}
\end{Lem}

\begin{proof} By (\ref{vj}) and (\ref{gamma-spec}),
it follows that the iterated Segre mapping $v^j$ is
of the form
\begin{equation}\Label{vj-decomp} v^{j}(t^1,\ldots,
t^{j})=(t^j,\nu^{j}(t^1,\ldots,t^j)),
\end{equation} where
\begin{equation}\Label{nuj}
\nu^j(t^1,\ldots,t^j)=Q(t^j,\bar
v^{j-1}(t^1,\ldots,t^{j-1})).
\end{equation} In particular, we obtain
\begin{multline}\Label{iterj+1} v^{j+1}(t^1,\ldots,
t^{j+1})=\\(t^{j+1},Q(t^{j+1},t^j,\bar
Q(t^j,t^{j-1},\nu^{j-1}(t^1,\ldots, t^{j-1})))).
\end{multline} The identity (\ref{hj})) is an
immediate consequence of (\ref{reality}). This
completes the proof of Lemma \ref{vj+1=vj-1}.
\end{proof}

We set, for $j\geq 1$,
\begin{equation}\Label{theta}
\Theta^j(t^1,\ldots,t^{j+1}):=(v^{j+1}(t^1,\ldots,
t^j,t^{j+1}+t^{j-1}),
\bar v^{j}(t^1,\ldots, t^j))
\end{equation} and, for $j\geq 2$,
\begin{equation}\Label{phi}
\Phi^j(t^1,\ldots,t^{j}):=(v^{j-1}(t^1,\ldots,
t^{j-1}),
\bar v^{j}(t^1,\ldots, t^j)).
\end{equation} Also, we define
\begin{equation}
\Theta^0(t^1):=(v^{1}(t^1),0),\quad
\Phi^1(t^1):=(0,
\bar v^{1}(t^1)).
\end{equation} Observe, by using (\ref{vj}), that,
for each $j\geq 1$, $\rk
\d v^{j+1} /\d t^{j+1}(0)=n$. Hence, it follows from
the  definition (\ref{theta}) of $\Theta^j$, that
\begin{equation}\Label{theta-rank}
\Rk\Theta^j=\Rk v^j+n, \quad j\geq 1,
\end{equation}
and similarly, 
\begin{equation}\Label{phi-rank}
\Rk\Phi^j=\Rk v^{j-1}+n, \quad j\geq 1.
\end{equation}
It follows from \S\ref{Prelim} that
$\Theta^j(\bC^{(j+1)n})\subset\mathcal M$, for $j\geq
1$, and $\Phi^j(\bC^{jn})\subset\mathcal M$, for
$j\geq 2$. By Lemma \ref{vj+1=vj-1}, we have, for
$j\geq 2$,
\begin{equation}\Label{phitheta}
\Theta^{j}(t^1,\ldots, t^j,0)=\Phi^j(t^1,\ldots,
t^j).
\end{equation} Also, by observing that
$\Rk\Theta^j=\Rk
\tilde\Theta^j$, where
\begin{multline}
\tilde\Theta^j(t^1,\ldots,t^{j+1}):=\Theta^j
(t^1,\ldots, t^{j+1}-t^{j-1})=\\ (v^{j+1}(t^1,\ldots,
t^j,t^{j+1}),
\bar v^{j}(t^1,\ldots, t^j)),
\end{multline} it is not difficult to see that
$\Rk\Theta^j=\Rk
\Phi^{j+1}$, for any $j\geq 0$.  A straightforward
calculation shows that, for any
$f\in\bC\dbl Z,\zeta\dbr$, any integer
$j\geq 0$, and for $l=1,\ldots,n$,
\begin{equation}\Label{push-1}
\begin{aligned}
\frac{\d}{\d t^{j+1}_l}(f\circ\Theta^j)(t^1,\ldots,
t^{j+1})= &\,((\tilde\mathcal L_l
f)\circ\Theta^j)(t^1,\ldots,t^{j+1}),
\\
\frac{\d}{\d t^{j+1}_l}(f\circ\Phi^{j+1})(t^1,\ldots,
t^{j+1})= &\,((\mathcal L_l
f)\circ\Phi^{j+1})(t^1,\ldots,t^{j+1}),
\end{aligned}
\end{equation}
where the $\mathcal L_j$ and $\tilde \mathcal L_j$
are the formal vector fields given by (\ref{Lj}).

\begin{Pro}\Label{orbit}Let
$\mathcal O$ denotes the formal CR orbit of
$\mathcal M$ (as defined in \S \ref{pfmain-3}). The
following hold.
\begin{enumerate}
\item [(i)]  
$\Theta^j(\bC^{(j+1)n})\subset\mathcal O$ and
$\Phi^{j+1}(\bC^{(j+1)n})\subset\mathcal O$, for
$j=0,1,\ldots$.
\item [(ii)] There exists $k_0$, $ 1\le k_0 \le
d+1$, such that
\begin{equation}\Label{rank-o}
\Rk
\Theta^{k_0+1}=\Rk \Phi^{k_0+1}=\dim \mathcal O. 
\end{equation} 
\end{enumerate} 
\end{Pro} 

\begin{proof}  By iterating (\ref{push-1}), we
obtain, for any multi-index
$\alpha\in\Bbb Z_+^n$,
\begin{equation}\Label{push-2}
\begin{aligned}
\left(\frac{\d}{\d
t^{j+1}}\right)^\alpha(f\circ\Theta^j)(t^1,\ldots,
t^{j+1})= &\,\big((\tilde\mathcal L^\alpha
f)\circ\Theta^j\big)(t^1,\ldots,t^{j+1}),
\\
\left(\frac{\d}{\d
t^{j+1}}\right)^\alpha(f\circ\Phi^{j+1})(t^1,\ldots,
t^{j+1})= &\,\big((\mathcal L^\alpha
f)\circ\Phi^{j+1}\big)(t^1,\ldots,t^{j+1}).
\end{aligned}
\end{equation} In particular, applying
(\ref{push-2}) with $j=0$, we conclude that, for
each multi-index $\alpha$,
\begin{equation}
\begin{aligned}
\left(\frac{\d}{\d
t^{1}}\right)^\alpha(f\circ\Theta^0)(0)=
&\,((\tilde\mathcal L^\alpha f)\circ\Theta^0)(0),
\\
\left(\frac{\d}{\d
t^{1}}\right)^\alpha(f\circ\Phi^{1})(0)=
&\,((\mathcal L^\alpha f)\circ\Phi^{1})(0).
\end{aligned}
\end{equation} Since $\mathcal L_l$, $\tilde
\mathcal L_l$ are tangent to
$\mathcal O$, we deduce that $$(\d /\d t^1)^\alpha
(f\circ\Theta^0)(0)=(\d /\d t^1)^\alpha
(f\circ\Phi^1)(0)=0$$ for all $f\in\mathcal
I(\mathcal O)$ and all multi-indices $\alpha$. Hence,
$f\circ\Theta^0=f\circ \Phi^1=0$, which proves that
$\Theta^0(\bC^n)\subset\mathcal O$ and
$\Phi^1(\bC^n)\subset\mathcal O$. 

Assume that
$\Theta^j(\bC^{(j+1)n})\subset\mathcal O$ and
$\Phi^{j+1}(\bC^{(j+1)n})\subset\mathcal O$, for all
$j=0,1,\ldots, j_0-1$. We shall prove it for
$j=j_0$. By applying (\ref{push-2}) with $j=j_0$, we
conclude, using also (\ref{phitheta}), that for
$f\in\mathcal I(\mathcal O)$,
\begin{equation}
\begin{aligned}
\left(\frac{\d}{\d
t^{{j_0+1}}}\right)^\alpha(f\circ\Theta^{j_0})
(t^1,\ldots,t^{j_0},0)= &\,((\tilde\mathcal L^\alpha
f)\circ\Theta^{j_0})(t^1,\ldots,t^{j_0},0)\\= &\,
((\tilde\mathcal L^\alpha
f)\circ\Phi^{j_0})(t^1,\ldots,t^{j_0})\\=&\,0,
\end{aligned}
\end{equation} where the last equality follows from
the inductive hypothesis and the fact that the
$\tilde\mathcal L_j$ are tangent to $\mathcal
I(\mathcal O)$. Hence,
$f\circ\Theta^{j_0}$ is independent of $t^{j_0+1}$.
Consequently,
$$ f\circ\Theta^{j_0}(t^1,\ldots,
t^{j_0+1})=f\circ\Theta^{j_0}(t^1,\ldots,
t^{j_0},0)=f\circ\Phi^{j_0}(t^1,\ldots, t^{j_0}),
$$ which is 0 by the inductive hypothesis. A similar
argument, using instead the identity
\begin{equation}\Label{phitheta-2}
\Phi^{j+1}(t^1,\ldots,
t^{j-1},t^j,t^{j+1})\big|_{t^{j+1}=t^{j-1}}=\tilde
\Theta^{j-1}(t^1,\ldots, t^j),
\end{equation}
 which follows from (\ref{hj}), shows that
$f\circ\Phi^{j_0+1}=0$ for
$f\in\mathcal I(\mathcal O)$. This completes the
proof of (i) of Proposition
\ref{orbit}.

To prove (ii) let $k_0$ be defined as in
\ref{strincr}).  By using
(\ref{incr}), (\ref{strincr}), (\ref{theta-rank}),
and (\ref{phi-rank}) we first observe that
$$
 \Rk \Theta^{j} < \Rk\Phi^{j},\ \  1 \le j \le k_0, \
\  {\rm and}\ \ \Rk \Theta^{j} = \Rk\Phi^{j}, \ \  j
\ge k_0 +1.
$$  In light of part (i), (\ref{phitheta}), and
(\ref{push-1}),
 the conclusion (\ref{rank-o}) of (ii) is an immediate
consequence of Proposition \ref{crucial} below with
$\Phi:=
\Phi^{k_0 +1}$ and $\Theta:=
\Theta^{k_0 +1}$.
\end{proof}

\begin{Pro}\Label{crucial} Let $\Sigma\subset\bC^m$
be a formal manifold and
$\Phi\colon(\bC^{p},0)\to (\bC^m,0)$,
$\Theta\colon(\bC^p\times\bC^q,0)\to (\bC^m,0)$
formal mappings such that the following hold.
\begin{enumerate}
\item [(i)]
$\Phi(x)=
\Theta(x,0),\quad x=(x_1,\ldots,x_p);
$
\item [(ii)]
$\Theta(\bC^p\times\bC^q)\subset\Sigma,
\ \ \Phi(\bC^p)\subset\Sigma$;
\item [(iii)]
$\Rk\Phi=\Rk\Theta$;
\item [(iv)]
There
are formal vector fields
$X_1,\ldots,X_k$, $Y_1,\ldots,$ $Y_k$ on $\bC^m$
tangent to $\Sigma$, formal vector fields $\hat
X_1,\ldots, \hat X_k$ on
$\bC^p$, and $\hat Y_1,\ldots,\hat Y_k$ on
$\bC^p\times\bC^q$ such that for every $f\in\bC\dbl
y\dbr$,
$y=(y_1,\ldots,y_m)$,
\begin{equation}\Label{push}
\begin{aligned}
\hat X_j(f\circ\Phi)=(X_jf)\circ \Phi &,\quad
j=1,\ldots, k,\\
\hat Y_j(f\circ\Theta)=(Y_jf)\circ \Theta &,\quad
j=1,\ldots, k,
\end{aligned}
\end{equation} and the vector space obtained by
 evaluating at
$0$ the elements of the
Lie algebra generated by
$X_1,\ldots,X_k,Y_1,\ldots, Y_k$  has dimension $\dim\Sigma$.
\end{enumerate}
 Then necessarily
$\Rk\Phi =\Rk\Theta= \dim\Sigma$.
\end{Pro}

\section{Proof of Proposition
$\ref{crucial}$}\Label{pfmain-5}

For the proof of Proposition \ref{crucial}, we first
make some preliminary reductions. Since
$\Sigma$ is a formal manifold through the origin, we
may choose the variables
$y_1,\ldots, y_m$ in $\bC^m$ so that the ideal
$\mathcal I(\Sigma)$ is generated by
$y_{s+1},\ldots, y_m$, where $s=\dim\Sigma$. Let us
write
$y=(y',y'')$, where $y'=(y_1,\ldots, y_s)$ and
$y''=(y_{s+1},\ldots, y_m)$. Similarly, we 
decompose the mappings $\Theta=(\Theta',\Theta'')$
and $\Phi=(\Phi',\Phi'')$, where e.g.\
$\Theta'=(\Theta_1,\ldots,\Theta_s)$ and
$\Theta_j=y_j\circ\Theta$. It is easy to see that
condition (ii) of the proposition
is equivalent to
$$
\Theta''(x,t)= 0,\quad\Phi''(x)= 0,\quad
x=(x_1,\ldots,x_p),\quad t=(t_1,\ldots, t_q).
$$ Moreover, a formal vector field $$X=\sum_{j=1}^s
a_j(y',y'')\frac{\d}{\d
  y_j}+\sum_{j=s+1}^m b_j(y',y'')\frac{\d}{\d y_j}$$
is tangent to
$\Sigma$ if and only if each $b_j(y',0)= 0$. Let us
write
$$ X'=\sum_{j=1}^s a_j(y',0)\frac{\d}{\d y_j}.
$$ Also, for $f\in\bC\dbl y',y''\dbr$, we write
$\tilde f(y'):=f(y',0)$. Then we have for $X$ as
above tangent to $\Sigma$,
$$ (Xf)\circ\Theta=
(X'\tilde f)\circ\Theta',\quad (Xf)\circ\Phi=
(X'\tilde f)\circ\Phi'.
$$
 Hence, if there exists a formal vector
field  $\hat X$ on $\bC^p$ such
that
$
 \hat
X(f\circ\Phi)=(Xf)\circ\Phi$ for all
$f\in\bC\dbl y\dbr,
$
 then  (since
$f\circ\Phi=
\tilde f\circ\Phi'$) it follows that
$
\hat X(\tilde f\circ\Phi')= (X'\tilde
f)\circ\Phi'$. 
 Similarly, if there exists a formal vector
field  $\hat Y$ on $\bC^p\times \bC^q$ such that $\hat
Y(f\circ\Theta)=(Xf)\circ\Theta$ for all
$f\in\bC\dbl y\dbr,$
then
$\hat Y(\tilde f\circ\Theta')= (X'\tilde
f)\circ\Theta'$.  Thus, by identifying
$\Sigma$ with
$\bC^s$, it suffices to prove Proposition
\ref{crucial} in the special case
$\Sigma =\bC ^m$.

For the proof of Proposition \ref{crucial} in the
case
$\Sigma = \bC^m$, we shall need some notation and
preliminary results. As before, we let
$\Bbb K_x$ denote the field of fractions of the ring
$\bC\dbl x\dbr$ of formal power series in the
variables
$x=(x_1,\ldots, x_n)$ and consider $\bC\dbl x\dbr$
as a subring of
$\Bbb K_x$. Recall that if
$F\colon(\bC^p,0)\to (\bC^m,0)$ is a formal mapping,
then $\Rk F$ denotes the rank of its Jacobian, i.e.\
the rank of the linear mapping $J_F\colon\Bbb
K^p_x\to \Bbb K^m_y$, where
$x=(x_1,\ldots, x_p)$ and $y=(y_1,\ldots, y_m)$,
defined by the $m\times p$-matrix
\begin{equation}\Label{cruc-5} J_F:=\pmatrix
F_{1,x_1} &\ldots&F_{1,x_p}\\ \vdots& &\vdots\\
F_{m,x_1} &\ldots&F_{m,x_p}\endpmatrix.
\end{equation} Recall that we denote by
$\varphi_F\colon\bC\dbl y\dbr\to \bC\dbl x\dbr$ the
homomorphism induced by $F$ and defined, for $f\in
\bC\dbl y\dbr$, by $\varphi_F(f):=f\circ F$. (Note
that e.g.\ the first equation of (\ref{push}) can
then be written $\hat X_j(\varphi_\Phi(f))=
\varphi_\Phi(X_jf)$.) 

The proof of Proposition \ref{crucial} in the case
$\Sigma =
\bC^m$ rests on the following four lemmas.
We shall use the notation and conventions
previously introduced. Recall in particular that
$\D(\bC\dbl x\dbr)$ denotes the set of all formal
vector fields in $x$.

\begin{Lem}\Label{lemma-7} Let $F\colon(\bC^p,0)\to
(\bC^m,0)$ be a formal mapping, $X\in\D(\bC\dbl
x\dbr)$,
$Y\in\D(\bC\dbl y\dbr)$ with
\begin{equation}
\begin{aligned} X &=\sum_{j=1}^p
a_j(x)\frac{\partial}{\partial x_j},\quad a_j\in
\bC\dbl x\dbr,\\ Y &=\sum_{j=1}^m
b_j(y)\frac{\partial}{\partial y_j},\quad b_j\in
\bC\dbl y\dbr,
\end{aligned}
\end{equation} and $c\in\bC\dbl x\dbr$,
$c\not= 0$.  Then the following are equivalent.
\begin{enumerate}
\item [(i)] 
$ X(\varphi_F(f))=
c\,\varphi_F(Yf)$ for all $f\in
\bC\dbl y\dbr$.
\item [(ii)] $X(\varphi_F(f_j)) =
c\,\varphi_F(Yf_j)$  
 for $f_j(y):= y_j$,  $\ j=1,\ldots, m.$
\item [(iii)]  $J_F\, a' = b\circ F,$ where
$$ a':=\pmatrix
a_1/c\\\vdots\\a_p/c\endpmatrix\in \Bbb K^p_x,\quad
b:=\pmatrix b_1\\\vdots\\b_m\endpmatrix\in \bC\dbl
y\dbr^m\subset\Bbb K^m_y.
$$
\end{enumerate} Moreover, any of the equivalent
conditions {\rm (i), (ii), (iii)} implies that the
formal vector field $Y$ is tangent to the ideal
$\ker\varphi_F$.
\end{Lem}

%%[{\it Remark.} Is the converse of the last
%%statement of Lemma 7
%%true? In other words, suppose that $Y$ is tangent to
%%$\ker\varphi_F$. Does this imply e.g.\ condition
%(i) above? I did
%%not spend much time on this since it is not (nor is
%the last
%%statement of Lemma 7 for that matter) relevant for
%the proof of
%%Proposition 1. However, it certainly seems relevant
%in the general
%%context of setting up an ideal based theory for
%Segre sets
%%etc.]\medskip

\begin{proof} The implication (i) $\implies$ (ii) is
clear. The opposite implication follows by an
inductive argument based on the facts that
$\varphi_F$ is a homomorphism, and $X$ and
$Y$ are derivations. For instance, we have, for
$f,g\in\bC\dbl y\dbr$,
$$ X(\varphi_F(fg))=
X(\varphi_F(f))\varphi_F(g)+\varphi_F(f)
X(\varphi_F(g)).
$$ A similar identity holds for
$\varphi_F(Y(fg))$. The details are left to the
reader.

The equivalence (ii) $\iff$ (iii) follows by simply
writing (ii) in terms of the components
$F_j=\varphi_F(y_j)$ of $F$.
The last statement of Lemma \ref{lemma-7} is a
direct consequence of (i).
\end{proof}

\begin{Lem}\Label{lemma-12} Let $F\colon(\bC^p,0)\to
(\bC^m,0)$ be a formal mapping, $X_1, X_2\in
\D(\bC\dbl x\dbr)$, $Y_1, Y_2\in
\D(\bC\dbl y\dbr)$ and suppose that there are
$c_1,c_2\in \bC\dbl x\dbr$, $c_j\not= 0$ for
$j=1,2$, such that
\begin{equation}\Label{cruc-13} X_j(\varphi_F(f))=
c_j\, \varphi_F(Y_jf),\quad \forall f\in\bC\dbl
y\dbr,\ j=1,2.
\end{equation} Then, there are $X'\in \D(\bC\dbl
x\dbr)$ and $c'\in\bC\dbl x\dbr$ with $c'\not= 0$,
such that
\begin{equation}\Label{cruc-14} X'(\varphi_F(f))=
c'\, \varphi_F([Y_1,Y_2]f),\quad
\forall f\in\bC\dbl y\dbr,
\end{equation} where $[\ ,\ ]$ denotes the
usual commutator of vector fields.
\end{Lem}

\begin{proof} A simple computation (which is left to
the reader) shows that (\ref{cruc-14}) holds with
$c'= (c_1c_2)^2$ and
\begin{equation}\Label{cruc-15} X'=
c_1c_2[X_1,X_2]+c_2(X_2c_1)X_1-c_1(X_1c_2)
X_2.\end{equation}
\end{proof}

The main idea in the proof of Proposition
\ref{crucial} is to reduce it to the following.

\begin{Lem}\Label {lemma-16} Let $F\colon(\bC^p,0)\to
(\bC^m,0)$ be a formal mapping, and $X_1,\ldots, X_k$
$\in \D(\bC\dbl x\dbr)$,
$Y_1,\ldots, Y_k\in \D(\bC\dbl y\dbr)$. Suppose that
there are
$c_1,\ldots, c_k\in \bC\dbl x\dbr$, $c_j\not = 0$ for
$j=1,\ldots k$, such that
\begin{equation}\Label{cruc-17} X_j(\varphi_F(f))=
c_j\, \varphi_F(Y_jf),\quad \forall f\in\bC\dbl
y\dbr,\ j=1,\ldots, k,
\end{equation} and such that the vector space obtained
by
 evaluating at
$0$ the elements of the
Lie algebra generated by
$Y_1,\ldots,Y_k$  has dimension $m$.
Then
$\Rk F=m$.
\end{Lem}

\begin{proof} Among $Y_1,\ldots, Y_k$ and all their
repeated commutators, we can pick out, by the
assumption in the lemma, $Y'_1,\ldots,
Y'_m\in\D(\bC\dbl y\dbr)$ such that
$Y'_1,\ldots, Y'_m$ evaluated at $0$ span
$T^{1,0}_0\bC^m$. By Lemma \ref{lemma-12}, we may
assume that there are $X'_1,\ldots,
X'_m\in\D(\bC\dbl x\dbr)$ and $c_1,\ldots, c_m\in
\bC\dbl x\dbr$,
$c_j\not = 0$ for $j=1,\ldots m$, such that
\begin{equation}\Label{cruc-18} X'_j(\varphi_F(f))=
c_j\, \varphi_F(Y'_jf),\quad \forall f\in\bC\dbl
y\dbr,\ j=1,\ldots, m.
\end{equation} Let us write
\begin{equation}\Label{cruc-19} Y'_j=\sum_{l=1}^m
b_{lj}(y)\frac{\partial}{\partial y_l},\quad
j=1,\ldots, m.\end{equation} By assumption the
vectors
\begin{equation}\Label{cruc-20} b_j(y):=\pmatrix
b_{1j}(y)\\\vdots\\ b_{mj}(y)\endpmatrix,\quad
j=1,\ldots m,\end{equation} evaluated at $0$ span
$\bC^m$. It is not difficult to see that this fact
implies that $\varphi_F(b_1),\ldots
\varphi_F(b_m)\in \bC\dbl x\dbr^m\subset \Bbb
K_x^m$ are linearly independent over $\Bbb K_x$ (and
hence span
$\Bbb K_x^m$) since they are linearly independent
over
$\bC$ when evaluated at
$0$. On the other hand, by Lemma \ref{lemma-7} and
(\ref{cruc-18}), these vectors are all in the image
of the Jacobian $J_F$, considered as a linear map from
$\Bbb
K_x^p$ to $\Bbb
K_x^m$.  Hence
$\Rk F=m$ and the lemma is proved.\end{proof}

The last lemma needed for the proof of Proposition
\ref{crucial} in the case $\Sigma=\bC^m$ is the
following.

\begin{Lem}\Label{lemma-21} Let $F\colon(\bC^p,0)\to
(\bC^m,0)$ be a formal mapping. Then the vector
subspace $J_F(\Bbb K_x^p)\subset \Bbb K^m_x$ is the
span over
$\Bbb K_x$ of the subset $J_F(\bC\dbl
x\dbr^p)\subset \Bbb K^m_x$.
\end{Lem}

\begin{proof} Let $a^1,\ldots, a^r\in \Bbb K^p_x$ be
such that $J_Fa^1,\ldots, J_Fa^r\in \Bbb K^m_x$ form
a basis for
$J_F(\Bbb K_x^p)$. By multiplying the $a_j$ by
suitable power series (clearing the denominators),
we obtain
$\tilde a^1,\ldots, \tilde a^r\in\bC\dbl x\dbr ^p$
and $J_F\tilde a^1,\ldots, J_F\tilde a^r\in \Bbb
K^m_x$  still form a basis for $J_F(\Bbb K_x^p)$
since
$J_F$ is linear over $\Bbb K_x$.
\end{proof}

\begin{proof}[Proof of Proposition $\ref{crucial}$
with $\Sigma = \bC^m$.] The main step in the proof
is to show that there are formal vector fields $\hat
Y'_1,\ldots,
\hat Y'_k$ on
$\bC^p$
and $c_1,\ldots, c_k\in
\bC\dbl x\dbr$, $c_j\not= 0$ for all $j$, such that
\begin{equation}\Label{cruc-22}
\hat Y_j'(f\circ\Phi)= c_j(Y_jf)\circ \Phi,\quad
\forall f\in
\bC\dbl y\dbr,\ j=1,\ldots k.\end{equation}
 To prove this, we write
$t=(t_1,\ldots, t_q)$ and denote by
$\eta\colon\bC\dbl x,t\dbr\to \bC\dbl x\dbr$ the
homomorphism defined by $\eta(g)(x):=g(x,0)\in
\bC\dbl x\dbr$ for $g\in\bC\dbl x,t\dbr$. By 
abuse of notation, we also denote by $\eta\colon
\bC\dbl x,t\dbr^m\to \bC\dbl x\dbr^m$ the mapping
given by applying $\eta$ to each component. By (i)
of Proposition  \ref{crucial}, for any
$a\in \bC\dbl x\dbr^p$, we have $\eta(J_\Theta a')=
J_\Phi a$, where
\begin{equation}\Label{cruc-23} a'(x,t):=\pmatrix
a_1(x)\\\vdots\\ a_p(x)\\0\\\vdots\\0\endpmatrix\in
\bC\dbl x,t\dbr^{p+q}.
\end{equation} Hence, we conclude that
$J_\Phi(\bC\dbl x\dbr^p)\subset
\eta(J_\Theta(\bC\dbl x,t\dbr^{p+q}))$, and therefore
by Lem\-ma~\ref{lemma-21}
\begin{equation}\Label{cruc-24}
J_\Phi(\Bbb K^p_x)\subset \text{\rm
span}_{\Bbb K_x}
\eta(J_\Theta(\bC\dbl x,t\dbr^{p+q})). \end{equation}
On the other hand, observe that if $v^1,\ldots,v^r
\in
\bC\dbl x,t\dbr^m$, with $\eta(v^1),\ldots,\eta(v^r)
 \in \bC\dbl
x\dbr^m$ linearly independent over $\Bbb K_x$,
then  
$v^1,\ldots,v^r$ are linearly independent over
$\Bbb K_{x,t}$.  Therefore we have
\begin{equation}\Label{cruc-29.5}
\dim\text{\rm span}_{\Bbb K_x} \eta(J_\Theta(\bC\dbl
x,t\dbr^{p+q}))\leq \dim\text{\rm span}_{\Bbb
K_{x,t}} (J_\Theta(\bC\dbl x,t\dbr^{p+q})).
\end{equation}  Hence, by using (\ref{cruc-24}),
(\ref{cruc-29.5}), Lemma \ref{lemma-21}, and the
assumption that
$\Rk\Theta=\Rk \Phi$, we obtain
\begin{multline}\Label{29.6}
\dim J_{\Phi}(\Bbb K_x^p)\leq 
\dim\text{\rm span}_{\Bbb K_x} \eta(J_\Theta(\bC\dbl
x,t\dbr^{p+q}))\leq \\\dim\text{\rm span}_{\Bbb
K_{x,t}} (J_\Theta(\bC\dbl x,t\dbr^{p+q}))=\dim
J_\Theta(\Bbb K_{x,t}^{p+q}))=\dim J_{\Phi}(\Bbb
K_x^p).
\end{multline} We therefore conclude again using
(\ref{cruc-24}) that
\begin{equation}\Label{cruc-30}
J_\Phi(\Bbb K_x^p)= \text{\rm
span}_{\Bbb K_x}
\eta(J_\Theta(\bC\dbl x,t\dbr^{p+q})).\end{equation}
Now we write
\begin{equation}\Label{cruc-31} Y_j=\sum_{l=1}^m
b_{l}^j(y)\frac{\partial}{\partial y_l},\quad
b_{l}^j\in\bC\dbl y\dbr,\ j=1,\ldots k,\end{equation}
and
\begin{equation}\Label{cruc-32} b^j(y):=\pmatrix
b_1^j(y)\\\vdots\\ b_{m}^j(y)\endpmatrix,\quad
j=1,\ldots k.
\end{equation} By Lemma
\ref{lemma-7} ((i)$\implies$(iii)) and the second
equation of (\ref{push}), we have
$b^j\circ\Theta\in J_\Theta(\bC\dbl x,t\dbr^{p+q})$,
and hence, by (\ref{cruc-30}) and (i) of Proposition
\ref{crucial},
$\eta(b^j\circ\Theta)=b^j\circ\Phi\in J_\Phi(\Bbb
K_x^p)$. By clearing denominators and again applying
Lemma
\ref{lemma-7} ((iii)$\implies$(i)), it follows
that there are
$\hat Y'_1,\ldots,
\hat Y'_k\in \D(\bC\dbl x\dbr)$ and $c_1,\ldots,
c_k\in \bC\dbl x\dbr$, $c_j\not= 0$, such that
(\ref{cruc-22}) holds. In view of the first
equation of (\ref{push}), the conclusion of
Proposition
\ref{crucial} for $\Sigma =\bC^m$  now follows by
applying Lemma
\ref{lemma-16} with $F=\Phi$.
\end{proof}

\section{Conclusion of the proof of Theorem
$\ref{formal-main}$; proof of Theorem
$\ref{formal-double}$}\Label{pfmain-6} 

The statements concerning the properties of
the iterated Segre mappings $v^j$ preceeding (i),
(ii) , and (iii) of Theorem
\ref{formal-main} have been proved in \S
\ref{pfmain-1} and \S\ref{pfmain-2}. Property (i)
has been proved in \S\ref{pfmain-3}, except for the
fact that $\mathcal O$ is CR, which will be proved
below.    We shall now prove statement (ii) of
Theorem
\ref{formal-main}. Let $\mathcal W$ be as in (ii) and
$g\in\mathcal I(\mathcal W)\subset
\bC\dbl Z\dbr$. We must show that $g\circ v^j=0$ for
each
$j\geq 1$. This follows immediately from Proposition
\ref{orbit} (i) since
$g$, considered as a power series in $(Z,\zeta)$
which is independent of $\zeta$, is also in
$\mathcal I(\mathcal O)$ and hence
$g\circ v^j=g\circ\Theta^{j-1}=0  $. 

To prove (iii), we use (\ref{phi-rank}) and 
(\ref{rank-o}) to obtain
\begin{equation} \Label {never}\Rk v^{k_0}=\Rk
\Phi^{k_0+1}-n=
\dim
\mathcal O-n.
\end{equation}
Recall from (i) and (ii) that $\dim \mathcal O = \dim
\frak g_M(0)$ and $\dim \mathcal W = N -e= \dim \frak
g_M(0) + d-N$. Hence (iii) follows from
(\ref{never}).

To complete the proof of Theorem \ref{formal-main},
it remains only to show that the formal manifold
$\mathcal O\subset\bC^N\times\bC^N$ is a CR
submanifold of $\bC^N$, i.e. that the ideal
$\mathcal I(\mathcal O)$ is real and
$\rk \d R/\d Z(0)$ equals the rank of $\d R/\d Z$
over the field of fractions $\Bbb K_{Z,\zeta}$. Here,
$R(Z,\zeta)=(\rho_1(Z,\zeta),\ldots,\rho_d(Z,\zeta),f_1(Z),\ldots,
f_e(Z))$ is a set of generators for the ideal
$\mathcal I(\mathcal O)$. To show reality, we let 
$g\in\mathcal I(\mathcal O)$. Recall that
$$
\sigma(g)(Z,\zeta):=\bar g(\zeta, Z).
$$ We must show that $\sigma(g)\in\mathcal
I(\mathcal O)$. Observe that, for each $j\geq 0$,
$$
\sigma(g)(\Phi^{j+1}(t^1,\ldots,
t^{j+1}))=\overline{g(\tilde\Theta^j(\bar
t^1,\ldots,\bar t^{j+1}))}.
$$ In view of Proposition \ref{orbit} (i), we have
$$ g(\tilde\Theta^j( t^1,\ldots,t^{j+1}))=0
$$ and, hence,
$$
\sigma(g)(\Phi^{j+1}(t^1,\ldots, t^{j+1}))=0.
$$ Since $\Rk\Phi^{k_0+1}=\dim\mathcal O$ by
(\ref{rank-o}), the conclusion $\sigma(g)\in\mathcal
I(\mathcal O)$ follows from Proposition
\ref{fact-1}. The fact that the formal real
submanifold $\mathcal O$ is CR is now a direct
consequence of the fact that $\mathcal M$ is generic
(and hence CR) and the fact that the remaining
generators $f_i$,
$i=1,\ldots,e$, are independent of
$\zeta$. This completes the proof of Theorem
\ref{formal-main}. \hfill$\square$

\begin{proof}[Proof of Theorem
$\ref{formal-double}$]  We first note that the
conclusion of the theorem is independent of the
choice of the Segre variety mapping $\gamma$. 
This can be seen by  using the notation of
\S\ref{pfmain-2} and observing
that
$v^{2k_0} =
\pi_1^{2k_0}\circ T^{2k_0}$, where
$T^{2k_0}$ is a parametrization of the formal
manifold $\frak S_{2k_0}$ . Hence it suffices to prove
Theorem
\ref{formal-double} using the special choice of
Segre variety mapping $\gamma$ introduced in
(\ref{gamma-spec}) in
\S
\ref{pfmain-4}. Let $\Sigma\subset
\bC^{2nk_0}$ be the formal manifold whose ideal is
generated by
$t^1$ and
$t^{2k_0-j}-t^{2+j}$, for $j=0,\ldots,k_0-2$. The
conclusion in Theorem \ref{formal-double} is now a
consequence of the definition of
$k_0$ in Theorem
\ref{formal-main} and of  Lemma 4.1.3 in \cite{MA}. 
\end{proof}

\end{document}